\documentclass{amsart}
\usepackage{amsmath,amssymb}
\usepackage{yhmath} 
\usepackage{mathrsfs}

\newcommand{\bdis}{\begin{displaymath}}
\newcommand{\edis}{\end{displaymath}}
\newcommand{\be}{\begin{equation}}
\newcommand{\ee}{\end{equation}}
\newcommand{\mbb}{\mathbb}
\newcommand{\mcal}{\mathcal}

\newcommand{\vp}{\varphi} 
\newcommand{\vth}{\vartheta}

\newcommand{\zf}{\zeta\left(\frac{1}{2}+it\right)}

\newcommand{\FR}{\frac{x^n+y^n}{z^n}} 
 
\newcommand{\FRm}{\frac{x^m+y^m}{z^m}}

\DeclareMathOperator{\re}{Re} 
\DeclareMathOperator{\im}{Im}

\theoremstyle{definition}

\theoremstyle{remark}
\newtheorem{remark}[]{Remark}

\newtheorem*{mydef11}{{\bf Theorem 1}}

\newtheorem*{mydef12}{{\bf Theorem 2}}

\newtheorem*{mydef13}{{\bf Theorem 3}}

\newtheorem*{mydef14}{{\bf Theorem 4}}

\newtheorem*{mydef15}{{\bf Theorem 5}}

\newtheorem*{mydef41}{{\bf Corollary 1}}

\newtheorem*{mydef51}{{\bf Lemma 1}}

\newtheorem*{mydef52}{{\bf Lemma 2}}

\newtheorem*{mydef81}{{\bf Property 1}}

\newtheorem*{mydef82}{{\bf Property 2}}

\newtheorem*{mydef91}{{\bf Formula 1}}

\newtheorem*{mydef92}{{\bf Formula 2}}

\newtheorem*{mydef93}{{\bf Formula 3}}

\newtheorem*{mydef94}{{\bf Formula 4}}

\numberwithin{equation}{section}

\begin{document}

\title[Jacob's ladders and wander from Hardy-Littlewood integral (1918) \dots]{Jacob's ladders and wander from Hardy-Littlewood integral (1918) to the $\zeta$-Platonic cube form of equivalent of the Fermat-Wiles theorem}

\author{Jan Moser}

\address{Department of Mathematical Analysis and Numerical Mathematics, Comenius University, Mlynska Dolina M105, 842 48 Bratislava, SLOVAKIA}

\email{jan.mozer@fmph.uniba.sk}

\keywords{Riemann zeta-function}

\begin{abstract}
In this paper we obtain $\zeta$-Platonic cube form of our almost linear formula. From this follows, for example, the $\zeta$-Platonic equivalent of the Fermat-Wiles theorem. The last result represents the point of contact between the Hardy-Littlewood integral and the Fermat-Wiles theorem and also with one of the Platonic solid - cube. 
\end{abstract}
\maketitle

\section{Introduction} 

\subsection{} 

In the paper \cite{12} we have defined the set 
\be \label{1.1} 
\{\mcal{A}_2(n)\}_{n=1}^{N(T)},\ n=1,\dots,n(T) 
\ee 
of the rectangular signals and we have obtained the following asymptotic formula 
\be \label{1.2} 
\sum_{n=1}^{N(T)}|\mcal{A}_2(n)|\sim (1-c)T,\ T\to\infty 
\ee 
on the base of Gramm's sequence $\{t_\nu\}$. 

\begin{remark}
	This formula represents the next sum variant of our almost linear formula, see \cite{10}, (3.4), (3.6), $r=1$ 
\be \label{1.3} 
\int_T^{\overset{1}{T}(T)}\left|\zf\right|^2{\rm d}t\sim (1-c)T,\ \overset{1}{T}(T)=\vp_1^{-1}(T),\ T\to\infty, 
\ee   
where the integral in (\ref{1.3}) represents new class of increments of the Hardy-Littlewood integral (1918), see \cite{4}. 
\end{remark} 

\subsection{} 

Now, let us remind that there are five regular convex polyhedra known as Platonic solids: 
\begin{itemize}
	\item[(i)] regular tetrahedron with four triangular faces, 
	\item[(ii)] regular hexahedron (cube) with six square faces, 
	\item[(iii)] regular octahedron with eight triangular faces, 
	\item[(iv)] regular icosahedron with twenty triangular faces, and 
	\item[(v)] regular dodecahedron with twelve pentagonal faces. 
\end{itemize} 
Here, of course, triangular face means equilateral triangle. In this connection see also the nice diagrams in \cite{2}, p. 226. 

In this paper we obtain new formulas generated by (\ref{1.2}) and containing areas of platonic cubes, infinitely small, in $\mbb{R}^3$.

\subsection{} 

We have obtained the following main formulas for the set 
\bdis 
 \left\{\mathscr{C}_1(n)\right\}_{n=1}^{N(T)}
\edis 
of Platonic cubes $ \mathscr{C}_1(n)$: $\zeta$-Platonic cube form: 
\begin{itemize}
	\item of almost linear formula (\ref{1.3}) 
	\be \label{1.4} 
	\sum_{T<t_n<\overset{1}{T}(T)}q(n)|\mathscr{C}_1(n)|\sim (1-c)T,\ T\to\infty, 
	\ee  
	\item of the first $\zeta$-functional\footnote{See \cite{11}, (4.5).} 
	\be \label{1.5} 
	\lim_{\tau\to\infty}\frac{1}{\tau}\left\{
	\sum_{t_n>\frac{x}{1-c}\tau}^{t_n<[\frac{x}{1-c}\tau]^1}q(n)|\mathscr{C}_1(n)|
	\right\} = x,\ x>0, 
	\ee 
	\item of equivalent of the Fermat-Wiles theorem 
	\be \label{1.6} 
	\lim_{\tau\to\infty}\frac{1}{\tau}\left\{
	\sum_{t_n>\FRm\frac{\tau}{1-c}}^{t_n<[\FRm\frac{\tau}{1-c}]^1}q(n)|\mathscr{C}_1(n)|
	\right\} \not= 1, 
	\ee   
	\item of asymptotic decomposition (or expansion) 
	\be \label{1.7} 
	\left|\mathscr{C}\left[\sqrt{\frac{\pi}{3}}\right]\right|\sim \sum_{n=1}^{N(T)}\frac{q(n)}{N(T)}|\mathscr{C}_1(n)|,\ T\to\infty, 
	\ee  
	of the area of Platonic cube with the edge $a=\sqrt{\frac{\pi}{3}}$ by means of areas of $\zeta$-Platonic cubes with infinitely small edges 
	\bdis 
	\sim \frac{\pi}{3\ln t_{\nu+n}},\ n=1,\dots,N(T),\ T\to\infty. 
	\edis 
\end{itemize} 

\begin{remark}
The $\zeta$-Platonic cube form of the equivalent of Fermat-Wiles theorem represents the point of contact between the Hardy-Littlewood integral (1918) and Fermat-Wiles theorem and also the Platonic hexahedron (cube). The first two notions are separated by the time interval more than 2200 years long from the last notion. 
\end{remark} 

\begin{remark}
The proof scheme that gives the results (\ref{1.4}) -- (\ref{1.7}) is adaptable also for the cases (i), (iii) -- (v) of the Platonic solids. 
\end{remark}

\section{Jacob's ladders: notions and basic geometrical properties}  

\subsection{}

In this paper we use the following notions of our works \cite{7} -- \cite{11}: 
\begin{itemize}
\item[{\tt (a)}] Jacob's ladder $\vp_1(T)$, 
\item[{\tt (b)}] direct iterations of Jacob's ladders 
\bdis 
\begin{split}
	& \vp_1^0(t)=t,\ \vp_1^1(t)=\vp_1(t),\ \vp_1^2(t)=\vp_1(\vp_1(t)),\dots , \\ 
	& \vp_1^k(t)=\vp_1(\vp_1^{k-1}(t))
\end{split}
\edis 
for every fixed natural number $k$, 
\item[{\tt (c)}] reverse iterations of Jacob's ladders 
\be \label{2.1}  
\begin{split}
	& \vp_1^{-1}(T)=\overset{1}{T},\ \vp_1^{-2}(T)=\vp_1^{-1}(\overset{1}{T})=\overset{2}{T},\dots, \\ 
	& \vp_1^{-r}(T)=\vp_1^{-1}(\overset{r-1}{T})=\overset{r}{T},\ r=1,\dots,k, 
\end{split} 
\ee   
where, for example, 
\be \label{2.2} 
\vp_1(\overset{r}{T})=\overset{r-1}{T}
\ee  
for every fixed $k\in\mbb{N}$ and every sufficiently big $T>0$. We also use the properties of the reverse iterations listed below.  
\be \label{2.3}
\overset{r}{T}-\overset{r-1}{T}\sim(1-c)\pi(\overset{r}{T});\ \pi(\overset{r}{T})\sim\frac{\overset{r}{T}}{\ln \overset{r}{T}},\ r=1,\dots,k,\ T\to\infty,  
\ee 
\be \label{2.4} 
\overset{0}{T}=T<\overset{1}{T}(T)<\overset{2}{T}(T)<\dots<\overset{k}{T}(T), 
\ee 
and 
\be \label{2.5} 
T\sim \overset{1}{T}\sim \overset{2}{T}\sim \dots\sim \overset{k}{T},\ T\to\infty.   
\ee  
\end{itemize} 

\begin{remark}
	The asymptotic behaviour of the points 
	\bdis 
	\{T,\overset{1}{T},\dots,\overset{k}{T}\}
	\edis  
	is as follows: at $T\to\infty$ these points recede unboundedly each from other and all together are receding to infinity. Hence, the set of these points behaves at $T\to\infty$ as one-dimensional Friedmann-Hubble expanding Universe. 
\end{remark}  

\subsection{} 

Let us remind that we have proved\footnote{See \cite{10}.} the existence of almost linear increments 
\be \label{2.6} 
\begin{split}
& \int_{\overset{r-1}{T}}^{\overset{r}{T}}\left|\zf\right|^2{\rm d}t\sim (1-c)\overset{r-1}{T}, \\ 
& r=1,\dots,k,\ T\to\infty,\ \overset{r}{T}=\overset{r}{T}(T)=\vp_1^{-r}(T)
\end{split} 
\ee 
for the Hardy-Littlewood integral (1918), \cite{4}: 
\be \label{2.7} 
J(T)=\int_0^T\left|\zf\right|^2{\rm d}t. 
\ee  

For completeness, we give here some basic geometrical properties related to Jacob's ladders. These are generated by the sequence 
\be \label{2.8} 
T\to \left\{\overset{r}{T}(T)\right\}_{r=1}^k
\ee 
of reverse iterations of the Jacob's ladders for every sufficiently big $T>0$ and every fixed $k\in\mbb{N}$. 

\begin{mydef81}
The sequence (\ref{2.8}) defines a partition of the segment $[T,\overset{k}{T}]$ as follows 
\be \label{2.9} 
|[T,\overset{k}{T}]|=\sum_{r=1}^k|[\overset{r-1}{T},\overset{r}{T}]|
\ee 
on the asymptotically equidistant parts 
\be \label{2.10} 
\begin{split}
& \overset{r}{T}-\overset{r-1}{T}\sim \overset{r+1}{T}-\overset{r}{T}, \\ 
& r=1,\dots,k-1,\ T\to\infty. 
\end{split}
\ee 
\end{mydef81} 

\begin{mydef82}
Simultaneously with the Property 1, the sequence (\ref{2.8}) defines the partition of the integral 
\be \label{2.11} 
\int_T^{\overset{k}{T}}\left|\zf\right|^2{\rm d}t
\ee 
into the parts 
\be \label{2.12} 
\int_T^{\overset{k}{T}}\left|\zf\right|^2{\rm d}t=\sum_{r=1}^k\int_{\overset{r-1}{T}}^{\overset{r}{T}}\left|\zf\right|^2{\rm d}t, 
\ee 
that are asymptotically equal 
\be \label{2.13} 
\int_{\overset{r-1}{T}}^{\overset{r}{T}}\left|\zf\right|^2{\rm d}t\sim \int_{\overset{r}{T}}^{\overset{r+1}{T}}\left|\zf\right|^2{\rm d}t,\ T\to\infty. 
\ee 
\end{mydef82} 

It is clear, that (\ref{2.10}) follows from (\ref{2.3}) and (\ref{2.5}) since 
\be \label{2.14} 
\overset{r}{T}-\overset{r-1}{T}\sim (1-c)\frac{\overset{r}{T}}{\ln \overset{r}{T}}\sim (1-c)\frac{T}{\ln T},\ r=1,\dots,k, 
\ee  
while our eq. (\ref{2.13}) follows from (\ref{2.6}) and (\ref{2.5}).

\section{The evolution of the Hardy-Littlewood integral from (1918) up to $\zeta$-Platonic hexahedron form of the Fermat-Wiles theorem} 

\subsection{} 

The origin of the Hardy-Littlewood integral lies in the asymptotic formula\footnote{See \cite{4}, pp. 122, 151 -- 156.} 
\be \label{3.1} 
\int_0^T\left|\zf\right|^2{\rm d}t\sim T\ln T,\ T\to\infty. 
\ee  
The Littlewood formula, see \cite{5}, $0<\delta$ being sufficiently small, 
\be \label{3.2} 
\int_0^T\left|\zf\right|^2{\rm d}t=T\ln T-(1+\ln 2\pi)T+\mcal{O}(T^{1/4+\delta})
\ee  
was the next step. And finally, Ingham proved the following formula, see \cite{6}, 
\be \label{3.3} 
\int_0^T\left|\zf\right|^2{\rm d}t=T\ln T-(1+\ln 2\pi - 2c)T+\mcal{O}(T^{1/2+\delta})
\ee 
that contains the complete set of main terms (here $c$ stands for the Euler's constant). 

\begin{remark}
We see that the first Ingham's contribution to the formula (\ref{3.2}) is the addition of the term $-2cT$ that is missing in previous Littlewood formula (\ref{3.2}). The second contribution of Ingham is quite essential improvement of the error term in Littlewood formula, namely the improvement is by 33.3\%. 
\end{remark} 

\begin{remark}
Let us remind also further considerable refinements of the error term, namely 
\be \label{3.4} 
\mcal{O}(T^{5/12+\delta}),\ 16.6\%;\quad \mcal{O}(T^{1/3+\delta}),\ 20\%
\ee 
due to Titchmarsh \cite{15} and Balasubramanian \cite{1}, correspondingly.  
\end{remark} 

We can conclude previous results into the following statement. 

\begin{mydef91}
\be \label{3.5} 
\int_0^T\left|\zf\right|^2{\rm d}t=T\ln T-(1+\ln 2\pi - 2c)T+R(T). 
\ee 
We will refer to (\ref{3.5}) as to Hardy-Littlewood-Ingham formula with the error term $R(T)$. 
\end{mydef91} 

Further, let us remind the fundamental Good's $\Omega$-theorem, see \cite{3}: 
\be \label{3.6} 
R(T)=\Omega(T^{1/4}),\ T\to\infty. 
\ee 
This statement gives us the exact terminal for the refinements of the type (\ref{3.4}). 

\begin{remark}
It is possible to write down 
\be \label{3.7} 
R(T)=\mcal{O}(T^{a+\delta}),\ \frac 14\leq a\leq \frac 13
\ee 
for the error term. In this form $a$ denotes the last valid exponent for which the estimate (\ref{3.7}) holds true. Especially, $a=\frac 14$ represents the Good's $\Omega$-theorem, and $a=\frac 13$ represents the main Barasubramanian's exponent. 
\end{remark} 

\begin{remark}
It is fully sufficient to use the exponent $\frac 13$ despite there is also the second Barasubramanian's exponent $\frac{27}{82}$ that is an improvement by approximately $1.18\%$ of $\frac 13$, see \cite{1}. 
\end{remark} 

\subsection{} 

Let us remind that we have introduced the Jacob's ladder 
\be \label{3.8} 
\vp_1(T)=\frac 12\vp(T)
\ee  
in the paper \cite{7}, comp. also \cite{8}, where the function $\vp(T)$ is an arbitrary solution to the non-linear integral equation 
\be \label{3.9} 
\int_0^{\mu[x(T)]}\left|\zf\right|^2e^{-\frac{2}{x(T)}t}{\rm d}t=\int_0^T\left|\zf\right|^2{\rm d}t, 
\ee  
where each admissible function $\mu(y)$ generates a solution 
\be \label{3.10} 
y=\vp(T;\mu)=\vp(T);\ \mu(y)\geq 7y\ln y. 
\ee  
We call the function $\vp_1(T)$ as Jacob's ladder since the analogy with the Jacob's dream in Chumash, Bereishis, 28:12. 

\begin{remark}
It is clear that the Jacob's ladder $\vp_1(T)$ represents new type of the special function generated by the modulus of the Riemann zeta-function on the critical line, i.e. 
\be \label{3.11} 
\left|\zf\right|\to \vp_1(t). 
\ee 
For the basic properties of this function see \cite{7} -- \cite{9}. 
\end{remark} 

Next, 83 years after the HLI formula (\ref{3.5}), we have proved the following result, see \cite{7}. 

\begin{mydef92}
The Hardy-Littlewood integral has, in addition to previously known HLI expression possessing an unbounded error term at $T\to\infty$,  the following  infinite number of new and almost exact representations 
\be \label{3.12} 
\begin{split}
& \int_0^T\left|\zf\right|^2{\rm d}t=\vp_1(T)\ln\{\vp_1(T)\}+ \\ 
& (c-\ln 2\pi)\vp_1(T)+c_0+\mcal{O}\left(\frac{\ln T}{T}\right),\ T\to\infty, 
\end{split}
\ee  
where $c$ is the Euler constant and $c_0$ is the constant from Titchmarsh-Kober-Atkinson formula.  
\end{mydef92} 

Further, we have obtained also the following (see \cite{7}, (6.2)): 
\begin{mydef93}
\be \label{3.13} 
T-\vp_1(T)\sim (1-c)\pi(T);\ \pi(T)\sim \frac{T}{\ln T},\ T\to\infty , 
\ee 
where the Jacob's ladder $\vp_1(T)$ can be seen as the complementary function to the function $(1-c)\pi(T)$ in the sense 
\be \label{3.14} 
\vp_1(T)+(1-c)\pi(T)\sim T,\ T\to\infty. 
\ee 
\end{mydef93} 

\begin{remark}
In the asymptotic formula (\ref{3.13}) we have used the following change\footnote{See also \cite{7}. (6.2).} 
\be \label{3.15} 
\frac{T}{\ln T}\to \pi(T),\ T\to\infty
\ee 
This is completely correct in the asymptotic region since the prime-number law. At the same time, we can give also a motivation for the change (\ref{3.15}). Namely, following\footnote{Comp. \cite{16}, p. 2.} 
\be \label{3.16} 
\zeta(s)=\exp\left\{s\int_0^{+\infty}\frac{\pi(x)}{x(x^s-1)}{\rm d}x\right\},\ \re\{s\}>1, 
\ee    
we can get that the function $\pi(x)$ is one element of the set of five functions generating the function $\zeta(s)$. 
\end{remark} 

\subsection{} 

Next, in the paper \cite{10}, we have proved the existence of the set of almost linear increments of the Hardy-Littlewood integral 
\be \label{3.17} 
J(T)=\int_0^T\left|\zf\right|^2{\rm d}t,\ T>0. 
\ee  
Namely, we have proved the following. 
\begin{mydef94}
For every fixed $k\in\mbb{N}$ and every sufficiently big $T>0$ it is true that\footnote{See (\ref{3.7}) about $a$ and $\delta$.}
\be \label{3.18} 
\begin{split}
& \int_{\overset{r-1}{T}}^{\overset{r}{T}}\left|\zf\right|^2{\rm d}t=(1-c)\overset{r-1}{T}+\mcal{O}(T^{a+\delta}), \\ 
& 1-c\approx 0.42,\ r=1,\dots,k,\ T\to\infty, 
\end{split}
\ee 
and\footnote{See (\ref{2.1}), (\ref{2.2}).} 
\be \label{3.19} 
\overset{r}{T}=\overset{r}{T}(T)=\vp_1^{-r}(T). 
\ee 
\end{mydef94} 

The final step of the proof of (\ref{3.18}) contains two main formulae only. Namely, by making use of the substitution $T\to \overset{r-1}{T}$ in (\ref{3.5}) and $T\to \overset{r}{T}$ in (\ref{3.12}) one obtains following formulae: 
\be \label{3.20} 
\int_{0}^{\overset{r-1}{T}}\left|\zf\right|^2{\rm d}t=\overset{r-1}{T}\ln\overset{r-1}{T}-(1+\ln2\pi-2c)\overset{r-1}{T}+\mcal{O}(T^{a+\delta}), 
\ee  
and 
\be \label{3.21} 
\begin{split}  
& \int_{0}^{\overset{r}{T}}\left|\zf\right|^2{\rm d}t=\overset{r-1}{T}\ln\overset{r-1}{T}-(c-\ln2\pi)\overset{r-1}{T}+c_0+\mcal{O}\left(\frac{\ln T}{T}\right), \\ 
& r=1,\dots,k,\ T\to\infty, 
\end{split}  
\ee 
correspondingly\footnote{See also (\ref{2.5}) and (\ref{2.2}).}. If we subtract (\ref{3.20}) from (\ref{3.21}), we obtain immediately (\ref{3.18}).  

\begin{remark}
It is clear that the formula (\ref{3.18}) is unreachable without the Jacob's ladder and its reverse iterations. 
\end{remark} 

\subsection{} 

Next step in the evolution of the Hardy-Littlewood integral represents the first $\zeta$-functional: 
\be \label{3.22} 
\lim_{\tau\to\infty}\frac{1}{\tau}\int_{\frac{x}{1-c}\tau}^{[\frac{x}{1-c}\tau]^1}\left|\zf\right|^2{\rm d}t=x,\ [G]^1=\vp_1^{-1}[G], 
\ee  
for every fixed $x>0$. This functional was obtained in the paper \cite{11}, (4.6). 

In the special case of the Fermat's rationals 
\bdis 
x\to \FR,\ x,y,z,n\in\mbb{N},\ n\geq 3
\edis 
we have obtained from (\ref{3.22}) also the first point of contact between almost linear increments (\ref{3.18}), $r=1$, and the Hardy-Littlewood integral and the Fermat-Wiles theorem. Namely, it is true: the $\zeta$-condition 
\be \label{3.23} 
\lim_{\tau\to\infty}\frac{1}{\tau}\int_{\FR\frac{\tau}{1-c}}^{[\FR\frac{\tau}{1-c}]^1}\left|\zf\right|^2{\rm d}t\not=1
\ee 
on the class of all Fermat's rationals represents the $\zeta$-equivalent of the Fermat-Wiles theorem. 

\subsection{} 

Present contribution to the evolution of the Hardy-Littlewood integral is expressed by the four $\zeta$-Platonic cube formulas (\ref{1.4}) -- (\ref{1.7}). 

\section{Next sum variant of almost linear formula} 

\subsection{} 

Let us remind the Riemann real-valued function 
\be \label{4.1} 
Z(t)=e^{i\vth(t)}\zf, 
\ee  
where\footnote{See \cite{13}, (35), (44), (62), comp. \cite{14}, p. 98.} 
\be \label{4.2} 
\begin{split}
& \vth(t)=-\frac t2\ln\pi+\im\ln\Gamma\left(\frac 14+i\frac t2\right)= \\ 
& \frac t2\ln\frac{t}{2\pi}-\frac t2-\frac{\pi}{8}+\mcal{O}\left(\frac 1t\right), 
\end{split}
\ee 
and, of course, 
\be \label{4.3} 
|Z(t)|=\left|\zf\right|. 
\ee 
In present paper we use also the Gramm' sequence 
\be \label{4.4} 
\begin{split}
& \{t_\nu\}_{\nu=1}^\infty:\ \vth(t_\nu)=\pi\nu, \\ 
& t_{\nu+1}-t_\nu=\frac{2\pi}{\ln t_\nu}+\frac{2\pi\ln2\pi}{\ln^2t_\nu}+\mcal{O}\left(\frac{1}{\ln^3t_\nu}\right), 
\end{split}
\ee 
where the last - Titchmarsh - formula can be found in \cite{14}, p. 102. In other words, we use the sequence of roots of the equations $\vth(t)=\pi\nu$. Next, we select from the sequence (\ref{4.4}) a local subsequence 
\be \label{4.5} 
\{t_n\}_{n=0}^{N+2},\ t_{\nu+n}\to t_n,\ N=N(T), 
\ee  
see \cite{12}, (3.8) and (3.9), where the letter $n$ denotes the local index, and 
\be \label{4.6}  
t_0\leq T<t_1<t_2<\dots<t_{N+1}<\overset{1}{T}(T)\leq t_{n+2}. 
\ee  
We use also the following formulas\footnote{See \cite{12}, (3.15) -- (3.18), (6.10).}: 
\be \label{4.7} 
Z^2(\tilde{t}(n))=\frac{1}{t_{n+1}-t_n}\int_{t_n}^{t_{n+1}}Z^2(t){\rm d}t,\ \tilde{t}(n)\in [t_n,t_{n+1}], 
\ee  
\be \label{4.8} 
\mcal{A}_2(n)=[t_n,t_{n+1}]\times[0,Z^2(\tilde{t}(n))],\ n=1,\dots,N(T), 
\ee  
\be \label{4.9} 
|\mcal{A}_2(n)|=Z^2(\tilde{t}(n))(t_{n+1}-t_n), 
\ee  
\be \label{4.10} 
\sum_{T<t_n<\overset{1}{T}(T)}|\mcal{A}_2(n)|\sim (1-c)T,\ T\to\infty, 
\ee  
\be \label{4.11} 
N(T)\sim \frac{1-c}{2\pi}T,\ T\to\infty, 
\ee 
correspondingly. 

\begin{remark}
One $\tilde{t}(n)\in [t_n,t_{n+1}]$ is chosen in the mean-value formula (\ref{4.7}). 
\end{remark} 

\subsection{} 

Now, we shall define the partition of $\mcal{A}_2(n)$ into the set of corresponding squares. 

For the segment $[t_n,t_{n+1}]$\footnote{See (\ref{4.8}).} we use the following partition 
\be \label{4.12} 
t_n=t_n^0<t_n^1<\dots<t_n^{r-1}<t_n^r<\dots<t_n^5<t_n^6=t_{n+1}, 
\ee  
where 
\be \label{4.13} 
t_n^r-t_n^{r-1}=\frac{1}{6}(t_{n+1}-t_n),\ r=1,\dots,6, 
\ee  
for every fixed $n=1,\dots,N(T)$. 

Now, we give the partition of the segment $[0,Z^2(\tilde{t}(n))]$\footnote{See (\ref{4.8}).} into the parts of length 
\bdis 
\frac 16(t_{n+1}-t_n). 
\edis  
\begin{itemize}
	\item[(A)] First, we shall consider such an $\bar{n}$ that 
	\be \label{4.14} 
	Z^2(\tilde{t}(\bar{n}))<\frac 16(t_{\bar{n}+1}-t_{\bar{n}}), 
	\ee  
	i.e. the case of \emph{bad} values of $n$. Since\footnote{See \cite{12}, (3.7), (6.3).} 
	\be \label{4.15} 
	t_{n+1}-t_n\sim \frac{2\pi}{\ln T},\ T\to\infty, 
	\ee  
	then\footnote{See (\ref{4.9}), (\ref{4.14}).} 
	\be \label{4.16} 
	|\mcal{A}_2(\bar{n})|<\frac 16(t_{n+1}-t_n)^2=\mcal{O}\left(\frac{1}{\ln^2T}\right),\ T\to\infty. 
	\ee 	
	Consequently, it is true\footnote{See (\ref{4.11}).} that 
	\be \label{4.17} 
	\begin{split}
	& \sum_{T<t_{\bar{n}}<\overset{1}{T}(T)}|\mcal{A}_2(\bar{n})|<\frac{A_1}{\ln^2T}\sum_{t_{\bar{n}}}1<\frac{A_1}{\ln^2T}N(T)< \\ 
	& A_2\frac{T}{\ln^2T}, 
	\end{split} 
	\ee 
	i.e. 
	\be \label{4.18} 
	\sum_{T<t_{\bar{n}}<\overset{1}{T}(T)}|\mcal{A}_2(\bar{n})|=\mcal{O}\left(\frac{T}{\ln^2T}\right)=o(T),\ T\to\infty. 
	\ee 
	\begin{remark}
	The meaning of the previous estimate is that the sum over all \emph{bad} values $\bar{n}$ gives negligible contribution into the sum (\ref{4.10}). Consequently, we can interpret the formula\footnote{See (\ref{4.10}).} 
	\be \label{4.19} 
	\sum_{T<t_{\bar{n}}<\overset{1}{T}(T)}|\mcal{A}_2(n)|\sim (1-c)T,\ T\to\infty, 
	\ee 
	as such one that holds true for the set of \emph{good} values $n$, i.e. for such $n$ for which it is true that 
	\be \label{4.20} 
	Z^2(\tilde{t}(n))\geq \frac 16(t_{n+1}-t_n),\ n=1,\dots,N_1(T)\sim N(T), 
	\ee 
	comp. (\ref{4.14}). After re-numbering the \emph{good} $n$s we can hold the old notation. 
	\end{remark} 
	\item[(B)] We define, in the case (\ref{4.20}), the partition of the segment $[0,Z^2(\tilde{t}(n))]$ by the set of the following points 
	\be \label{4.21} 
	y_n^0=0<y_n^1<\dots<y_n^{p-1}<y_n^p<\dots<y_n^{q-1}<y_n^q<y_n^{q+1}, 
	\ee 
	where 
	\be \label{4.22} 
	y_n^p-y_n^{p-1}=\frac{1}{6}(t_{n+1}-t_n),\ p=1,\dots,q(n), 
	\ee  
	\be \label{4.23} 
	y_n^q\leq Z^2(\tilde{t}(n))<y_n^{q+1}, 
	\ee  
	and 
	\be \label{4.24} 
	\frac{Z^2(\tilde{t}(n))}{\frac 16(t_{n+1}-t_n)}=q(n)+R(n),\ q(n)\in\mbb{N},\ R(n)\in[0,1]. 
	\ee  
	Consequently\footnote{See (\ref{4.24}), (\ref{4.15}).}, it is true 
	\be \label{4.25} 
	Z^2(\tilde{t}(n))=\frac{1}{6}(t_{n+1}-t_n)q(n)+\mcal{O}\left(\frac{1}{\ln T}\right),\ T\to\infty 
	\ee 
	for every fixed $n=1,\dots,N(T)$, see (\ref{4.20}). 
\end{itemize} 

\subsection{} 

Now, we define the following set of squares 
\be \label{4.26} 
\begin{split}
& \mcal{A}_3^{r,p}(n)=[t_n^{r-1},t_n^r]\times[y_n^{p-1},y_n^p], \\ 
& r=1,\dots,6,\ p=1,\dots,q(n)
\end{split}
\ee 
as the elements of the partition of $\mcal{A}_2(n)$, where\footnote{See (\ref{4.13}), (\ref{4.22}), (\ref{4.26}).} 
\be \label{4.27} 
|\mcal{A}_3^{r,p}(n)|=\frac{1}{36}(t_{n+1}-t_n)^2, 
\ee  
and 
\be \label{4.28} 
\mcal{A}_3(n)=\bigcup_{r=1}^6\bigcup_{p=1}^{q}\mcal{A}_3^{r,p}(n), 
\ee  
\be \label{4.29} 
|\mcal{A}_3(n)|=\sum_{r=1}^6\sum_{p=1}^{q}|\mcal{A}_3^{r,p}(n)|=\frac{1}{6}(t_{n+1}-t_n)^2q 
\ee  
for every fixed $n=1,\dots,N(T)$, see Remark 13, (\ref{4.20}). 

Next, it follows from (\ref{4.25})\footnote{See (\ref{4.15}).} 
\be \label{4.30} 
Z^2(\tilde{t}(n))(t_{n+1}-t_n)=\frac 16(t_{n+1}-t_n)^2q(n)+\mcal{O}\left(\frac{1}{\ln^2T}\right). 
\ee 
Consequently, the following lemma holds true\footnote{See (\ref{4.9}), (\ref{4.29}), (\ref{4.30}).}. 

\begin{mydef51}
\be \label{4.31} 
\begin{split}
& |\mcal{A}_2(n)|=|\mcal{A}_3(n)|+\mcal{O}\left(\frac{1}{\ln^2T}\right), \\ 
& n=1,\dots,N(T),\ T\to\infty. 
\end{split}
\ee 
\end{mydef51}  

And finally, we have the following result. 

\begin{mydef11}
\be \label{4.32} 
\sum_{T<t_n<\overset{1}{T}(T)}|\mcal{A}_3(n)|\sim (1-c)T,\ T\to\infty. 
\ee 
\end{mydef11} 

\begin{remark}
The formula (\ref{4.32}) expresses the next sum variant of our almost linear formula, see (\ref{2.6}), $r=1$. 
\end{remark} 

\section{The $\zeta$-Platonic cube form of equivalent of the Fermat-Wiles theorem} 

\subsection{} 

Let us begin with the next remark. 

\begin{remark}
The subsection 1.2 of this paper contains the set of five existence theorems. For example, by (ii) it is true: the Platonic cube therein is for every $a>0$ as a length of side for every six congruent and generating squares. Simultaneously, $a$ is the edge of the cube.  
\end{remark} 

Further, we introduce the Cartesian coordinate system $O_{tyz}$ in $\mbb{R}^3$. We use the positive octant 
\be \label{5.1} 
O^1_{tyz}:\ t\geq 0, y\geq 0, z\geq 0, 
\ee  
where 
\be \label{5.2} 
\mcal{A}_2(n)\subset O^1_{txy},\ n=1,\dots,N(T). 
\ee  
Since by (ii) there is the cube $\mathscr{C}[\alpha]$ with the edge\footnote{See (\ref{4.13}), (\ref{4.22}).}
\bdis 
\alpha=\frac{1}{6}(t_{n+1}-t_n), 
\edis  
we assign to the set 
\be \label{5.3} 
\left\{\bigcup_{r=1}^6\mcal{A}_3^{r,p}(n)\right\}\subset O^1_{ty}\subset\mbb{R}^2
\ee 
for every fixed $p$ and $n$  the Platonic cube (its surface) 
\be \label{5.4} 
\mathscr{C}[\alpha]=\mathscr{C}_p(n)\subset O^1_{tyz}\subset\mbb{R}^3, 
\ee  
with the property\footnote{Conservation of measures - in our case areas.}
\be \label{5.5} 
|\mathscr{C}_p(n)|=\sum_{r=1}^6|\mcal{A}_3^{r,p}(n)|. 
\ee 

Consequently, there is, for our case, the complete set 
\be \label{5.6} 
\left\{\mathscr{C}_p(n)\right\}_{p=1,n=1}^{q(n),N(T)}
\ee 
of $\zeta$-Platonic cubes in the set $O^1_{tyz}$ in $\mbb{R}^3$. 

The location of elements (\ref{5.6}) over 
\bdis 
\mcal{A}_2(n)\ (\supset \mcal{A}_3(n)), 
\edis  
that is for every fixed $n$, is controlled as follows. 

\begin{remark}
Namely, the location of elements 
\be \label{5.7} 
\left\{\mathscr{C}_p(n)\right\}_{p=1}^{q(n)}\subset \mcal{A}_n(n)\times [0,+\infty)\subset O^1_{tyz}
\ee  
in corresponding subset fulfils the conditions: 
\begin{itemize}
	\item[(a)] 
	\be \label{5.8} 
	\begin{split}
	& \mathscr{C}_p(n)\cup \mathscr{C}_{p'}(n)=\emptyset,\ p\not= p', \\ 
	& 1\leq p,p'\leq q(n),\ n=1,\dots,N(T), 
	\end{split}
	\ee  
	\item[(b)] for the orthogonal projection of every element $\mathscr{C}_p(n)$ in the direction $O^1_z$, i.e.  
	\bdis 
	[\mathscr{C}_p(n)]_\perp\to O^1_{ty}, 
	\edis 
	it is true that 
	\be \label{5.9} 
	[\mathscr{C}_p(n)]_\perp\subset \mcal{A}_2(n). 
	\ee 
\end{itemize}
\end{remark} 

\subsection{} 

Since the set 
\be \label{5.10} 
\left\{\mathscr{C}_p(n)\right\}_{p=1}^{q(n)}
\ee 
contains for every fixed $n=1,\dots,N(T)$ congruent elements only, we give the new results for the simple case 
\be \label{5.11} 
\sum_{p=1}^{q(n)}|\mathscr{C}_p(n)|=q(n)|\mathscr{C}_1(n)|=|\mcal{A}_3(n)|, 
\ee 
for example, that is sufficient. In this case we have obtained the following results. 

First, it is true\footnote{See (\ref{4.32}), (\ref{5.11}).}: 

\begin{mydef12}
\be \label{5.12} 
\sum_{T<t_n<\overset{1}{T}(T)}q(n)|\mathscr{C}_1(n)|\sim(1-c)T,\ T\to\infty. 
\ee 
\end{mydef12}  

\begin{remark}
The formula (\ref{5.12}) expresses the first $\zeta$-Platonic cube form of or almost linear formula, see (\ref{2.6}), $r=1$. 
\end{remark} 

\begin{remark}
From the geometric point of view ($|\mathscr{C}_1(n)|$ is area) more exact form of (\ref{5.12}) is as follows: 
\bdis 
\begin{split}
& \sum_{T<t_n<\overset{1}{T}(T)}q(n)|\mathscr{C}_1(n)|\sim(1-c)|\mcal{T}|,\ T\to\infty, \\ 
& \mcal{T}=[0,T]\times [0,1];\ |\mcal{T}|=T. 
\end{split} 
\edis  
We shall, however, use the form (\ref{5.12}). 
\end{remark} 

Now, if we use the substitution 
\be \label{5.13} 
T=\frac{x}{1-c}\tau,\ \{T\to\infty\} \ \Leftrightarrow \ \{\tau\to\infty\}
\ee 
in (\ref{5.12}), then we obtain the following result. 

\begin{mydef13}
\be \label{5.14} 
\lim_{\tau\to\infty}\frac{1}{\tau}\left\{
\sum_{t_n>\frac{x}{1-c}\tau}^{t_n<[\frac{x}{1-c}\tau]^1}q(n)|\mathscr{C}_1(n)|
\right\}=x 
\ee  
for every fixed $x>0$. 
\end{mydef13} 

\begin{remark}
The formula (\ref{5.14}) represents the $\zeta$-Platonic cube form of the first $\zeta$-functional\footnote{See (\ref{3.22}) and \cite{11}, (4.6).}.
\end{remark} 

In the special case of of the Fermat's rationals 
\bdis 
x\to\FRm,\ x,y,z,m\in\mbb{N},\ m\geq 3, 
\edis  
we have obtained the following two results. 

\begin{mydef41}
\be \label{5.15} 
\lim_{\tau\to\infty}\frac{1}{\tau}\left\{
\sum_{t_n>\FRm\frac{\tau}{1-c}}^{t_n<[\FRm\frac{\tau}{1-c}]^1}q(n)|\mathscr{C}_1(n)|
\right\}=\FRm
\ee 
for every fixed Fermat's rational. 
\end{mydef41} 

\begin{mydef14}
The $\zeta$-condition 
\be \label{5.16} 
\lim_{\tau\to\infty}\frac{1}{\tau}\left\{
\sum_{t_n>\FRm\frac{\tau}{1-c}}^{t_n<[\FRm\frac{\tau}{1-c}]^1}q(n)|\mathscr{C}_1(n)|
\right\}\not=1
\ee 
on the set of all Fermat's rationals expresses the $\zeta$-Platonic cube form of equivalent of the Fermat-Wiles theorem based on the Gramm' sequence $\{t_\nu\}$. 
\end{mydef14} 

Next, if we rewrite the eq. (\ref{5.12}) in the alternative form\footnote{Comp. \cite{12}, (6.11).} 
\be \label{5.17} 
\sum_{n=1}^{N(T)}q(n)|\mathscr{C}_1(n)|\sim (1-c)T,\ T\to\infty, 
\ee  
and remind that\footnote{See (\ref{4.11}).} 
\bdis 
N(T)\sim\frac{1-c}{2\pi}T,\ T\to\infty, 
\edis  
we obtain the following lemma. 

\begin{mydef52}
\be \label{5.18} 
\sum_{n=1}^{N(T)}\frac{q(n)}{N(T)}|\mathscr{C}_1(n)|\sim 2\pi,\ T\to\infty. 
\ee 
\end{mydef52} 

Now, we can introduce the Platonic cube 
\bdis 
\mathscr{C}\left[\sqrt{\frac{\pi}{3}}\right]
\edis 
with the edge $a=\frac{\pi}{3}\approx 1.02$ where, of course, 
\be \label{5.19} 
\left|\mathscr{C}\left[\sqrt{\frac{\pi}{3}}\right]\right|=2\pi. 
\ee 
Consequently, we obtain the following theorem. 

\begin{mydef15}
The asymptotic formula 
\be \label{5.20} 
\left|\mathscr{C}\left[\sqrt{\frac{\pi}{3}}\right]\right|\sim \sum_{n=1}^{N(T)}\frac{q(n)}{N(T)}|\mathscr{C}_1(n)|,\ T\to\infty 
\ee 
expresses decomposition (or expansion) of the area of Platonic cube with the edge $\sqrt{\frac{\pi}{3}}$ by means of the areas of $\zeta$-Platonic cubes with infinitely small edges\footnote{Comp. (\ref{4.4}), (\ref{4.5}).}
\bdis 
\sim \frac{\pi}{3\ln t_{\nu+n}},\ n=1,\dots,N(T),\ T\to\infty , 
\edis  
and also by the set 
\bdis 
\left\{\frac{q(n)}{N(T)}\right\}_{n=1}^{N(T)} 
\edis 
of rational numbers. 
\end{mydef15}

I would like to thank Michal Demetrian for his moral support of my study of Jacob's ladders.

\end{document}